\documentclass[12pt,a4paper,draft]{article}
\usepackage{amsmath}
\usepackage{amsfonts}
\usepackage{amssymb}
\makeatletter\@addtoreset{equation}{section}\makeatother

\newtheorem{theorem}{Theorem}[section]

\newtheorem{lemma}{Lemma}
\def\dd{\partial}
\def\g{{\frak g}}
\def\C{{\mathbb C}}

\newtheorem{definition}{Definition}

\begin{document}
\title{{\bf\large
QUANTIZATION OF GEOMETRIC CLASSICAL {\LARGE r}-MATRICES}}
\author{ {\Large Pavel Etingof and Alexandre Soloviev}\\
Harvard University\\
Department of Mathematics\\
Cambridge, MA 02138, USA\\
and\\
MIT \\
Department of Mathematics\\
Cambridge, MA 02139,
      USA
\date{March 1, 1999}
}
\maketitle

In this note we define geometric classical r-matrices and quantum 
R-matrices, and show how any geometric classical r-matrix can be quantized 
to a geometric quantum R-matrix. This is one of the simplest 
nontrivial examples of quantization of solutions of the classical Yang-Baxter 
equation, which can be explicitly computed. The idea of the above 
quantization was inspired by the results in \cite{ESS}. 
We note that a construction similar to ours was obtained in \cite{KLM}.  

\section{Geometric classical r-matrices and quantum R-matrices}

Let $X$ be a smooth, affine algebraic variety over $\mathbb C$. 

\begin{definition} A geometric classical r-matrix
on $X$ is a derivation \linebreak 
$r: \mathbb C[X\times X]\to \mathbb C[X\times X]$
(i.e. a vector field on $X\times X$), 
which satisfies the classical Yang-Baxter equation
\begin{equation}
\label{cybe}
[r_{12},r_{13}]+[r_{12},r_{23}]+[r_{13},r_{23}]=0\ {\text in} \ \mathbb
C[X\times X
\times X]
\end{equation}
and the unitarity condition 
\begin{equation}
\label{skew}
r+r_{21}=0 \ {\text in} \ \mathbb C[X\times X].
\end{equation}
\end{definition}

{\bf Example 1.} Let $X$ be any variety as above, and $v$ a vector field 
on $X$. Define $r^v(x,y)=(v(x),-v(y))$. Then $r$ is a geometric 
classical r-matrix. We call it a permutation r-matrix, since it corresponds 
to an ``infinitesimal permutation'' of $X$ given by $v$.  

{\bf Example 2.} Let $X$ be a finite dimensional algebra over $\mathbb C$
(e.g. a matrix algebra), and the vector field $r$ be given by 
$r_c(x,y)=(xcy,-ycx)$, where $c\in X$. It can be checked that 
$r_c$ is a geometric classical r-matrix.  

\begin{definition} A formal diffeomorphism 
$g$ of a smooth affine variety $Y$ 
is an algebra homomorphism $g:\C[Y]\to \C[Y][[\hbar]]$
such that $g=1+O(\hbar)$. 
\end{definition}

In particular, if $v$ is a vector field on $Y$, then one can define 
a formal diffeomorphism $g=e^{\hbar v}$ of $Y$ by 
$(gF)(x)=\sum_{m\ge 0}\frac{\hbar^mv^m}{m!}F(x)$. 
The last expression can be written as $F(e^{\hbar v}x)$, 
where $e^{\hbar v}x$ is a regular map from the formal disk to $Y$. 

\begin{definition} A geometric quantum R-matrix
on $X$ is a formal diffeomorphism of $X\times X$, 
which satisfies the quantum Yang-Baxter equation
\begin{equation}
\label{QYBE}
R_{12}R_{13}R_{23}=R_{23}R_{13}R_{12}
\end{equation}
and the unitarity condition 
\begin{equation}
RR_{21}=1.
\end{equation}
\end{definition}

This definition is a modification of Drinfeld's definition of a (unitary) 
set-theoretical solution of the quantum Yang-Baxter equation (see \cite{Dr}), 
in the case when $X$ is an algebraic variety. The term ``geometric'' is used 
because the map $R:\C[X^2]\to \C[X^2][[\hbar]]$ is not an arbitrary linear 
map, but a map of geometric origin, i.e. coming from a formal diffeomorphism 
of $X^2$.   

{\bf Example 3.} Let $X$ be as in Example 1. 
For any formal diffeomorphism $g$, define $R^g(x,y)=(g(x),g^{-1}(y))$.
This is a geometric 
quantum R-matrix. We call it a permutation R-matrix, since it corresponds 
to a ``formal permutation'' of $X$ given by $g$.  

{\bf Example 4.} Let $X$ be a finite dimensional algebra over $\mathbb C$, 
and the formal diffeomorphism $R$ be given by 
$R_c(x,y)=(x(1+\hbar cy),y(1+\hbar cx+\hbar^2cxcy)^{-1})$, 
where $c\in X$. It was checked in \cite{ESS} (see formula (A5)) that 
$R$ is a geometric quantum R-matrix.  

\vskip .05in

Suppose that $R$ is a geometric quantum R-matrix on $X$, and 
its $\hbar$-expansion looks like $R=1+\hbar r+O(\hbar^2)$. 
Then it is easy to check that $r$ is a geometric classical r-matrix. 

\begin{definition} $R$ is said to be a quantization of $r$. 
\end{definition}

{\bf Example 5.} It is easy to see that 
$R^g$ is a quantization of $r^v$ if $g=e^{\hbar v}$, and 
$R_c$ is a quantization of $r_c$.

Our main result is the following quantization theorem. 

\begin{theorem}
Any geometric classical r-matrix admits a quantization. 
\end{theorem} 

{\bf Remark.} It was proved in \cite{EK} that any classical r-matrix can be 
quantized. In the unitary case, this was proved earlier by Drinfeld. 
However, these results don't automatically guarantee that if 
the r-matrix $r$ is geometric then it has a quantization which is also 
geometric. So the main theorem does not obviously follow from 
the general quantization results. Also, the main theorem has 
the advantage that its proof gives an easy way to compute the 
quantization.  

\section{Proof of the main theorem}

{\bf 2.1. The Lie algebra with a bijective 1-cocycle associated to a geometric 
classical r-matrix.}

Let $r$ be a geometric classical r-matrix on $X$. Then $r$ is a vector 
field on $X^2$, so it is an element of 
$\text{Vect}X  \otimes   \mathbb C[X]  \oplus  \mathbb
C[X]  \otimes   
\text{Vect}(X)$, where $\text{Vect}(X)$ is the Lie algebra of vector fields
on $X$. \linebreak Consider the space 
$\g=\{(1\otimes f)(r)|f\in (\text{Vect}(X))^*\oplus (\mathbb C[X])^*\}$.
It follows from (\ref{cybe}) -- (\ref{skew}) that $\g$ is a 
finite dimensional Lie subalgebra
in the Lie algebra 
$\text{Vect}(X)\ltimes \mathbb C[X]$ 
of differential operators of order 
$\le 1$ on $X$. Moreover, $\g=\g_+\oplus \g_-$ as a
vector
space, where
$\g_+=\{(1\otimes f)(r)|f\in (\mathbb C[X])^*\}$,
$\g_-=\{(f\otimes 1)(r)|f\in (\text{Vect}[X])^*\}$ are Lie subalgebras, and 
$r\in \g_+\otimes \g_-\oplus \g_-\otimes \g_+$. 
Since
$[\g_+,\g_-]\subset \g_-$, the space $V=\g_-$ has a $\g_+$-module structure. 
We introduce a bijective map $\phi_r: V^*\to \g_+$ by the
formula $\phi_r(f)=(1\otimes f)(r)$ (here 
$f$ is extended by $0$ to $\g_+$). Denote $\pi=\phi_r^{-1}:\g_+\to V^*$.
 
\begin{lemma}
$\pi:\g_+\to V^*$ is a bijective 1-cocycle. That is, for any $a,b\in \g_+$,
$\pi([a,b])=a*\pi(b)-b*\pi(a)$, where $*$ denotes the $\g_+$ action on
$V^*$.  
\end{lemma}

{\it Proof of the Lemma.}

Let $f, g\in V^*$, $x\in \g_+^*$, then 
\begin{multline*}
[\phi_r(f),\phi_r(g)](x)=(x\otimes f\otimes g)([r_{12},r_{13}])\\
=-(x\otimes f \otimes g)([r_{12},r_{23}]+[r_{13},r_{23}])\\
=-f([(x\otimes 1)(r),(1\otimes g)(r)])+g([(x\otimes 1)(r),(1\otimes
f)(r)])\\
=-(\phi_r(g)*f)(x)+(\phi_r(f)*g)(x),
\end{multline*}
which proves the Lemma.
$\square$

{\bf 2.2. Exponentiation of the bijective cocycle.} 

Recall some basic facts about formal groups. 
Let $L$ be any Lie algebra over $\C$. We denote by 
$E(L)$ the 
group of formal expressions of the form $e^{\hbar b}$, where
$b\in L[[\hbar]]$, which are multiplied by the Campbell-Hausdorff 
formula. This is the group of $\C[[\hbar]]$-rational points of the 
formal group associated to the Lie algebra $L$. 

It is clear that $E$ is a functor from the category of Lie algebras 
to the category of groups. That is, to any homomorphism 
$\phi:L\to L'$ of Lie algebras there corresponds a homomorphism 
of groups $E(\phi):E(L)\to E(L')$. 

Let $\psi:L\to Vect(X)$ be a homomorphism of Lie algebras. Then the formal
group $E(L)$ acts on $X$ {\it on the right} by formal diffeomeorphisms, 
via $e^{hb}\to e^{h\psi(b)}$. We stress that this is a right, and not
left, action, i.e. the above assignment is an antihomomorphism, rather
than homomorphism, of groups. 

Now let us come back to the proof of 
the theorem. Define a linear map $\bar\pi:\g_+\to \g_+\ltimes V^*$ by 
$\bar\pi(a)=(a,\pi(a))$. Lemma 1 states that $\bar\pi$ is
a homomorphism of Lie algebras. 

Let $G_+=E(\g_+)$. Since $\g_+$ is a Lie algebra 
of vector fields on $X$, the group $G_+$ acts on $X$ {\it on the right} by
formal diffeomorphisms. We denote this action by $\rho$,
i.e. $F(\rho(e^{hb})x)=e^{hb}F(x)$ for $b\in \g_+,\ F\in \C[X]$. 

Also, it is clear that $G_+$ acts naturally on
$V^*[[\hbar]]$. 
Consider the group $G_+\ltimes V^*[[\hbar]]$ 
with multiplication 
$(a,b)(a',b')=(aa',b + a*b')$. It is easy to see that this group 
is naturally isomorphic to  
$E(\g_+\ltimes V^*)$ (here $\ltimes$ is the semidirect product). Therefore,   
the Lie algebra homomorphism $\bar\pi$ can be lifted to a group
homomorphism $\bar\Pi:G_+\to G_+\ltimes V^*[[\hbar]]$. 
Let $p:G_+\ltimes V^*[[\hbar]]\to V^*[[\hbar]]$ 
be the projection map. The bijective map
$\Pi:G_+\to V^*[[\hbar]]$
defined as a composition $\Pi=\hbar^{-1}p\bar\Pi$ satisfies 
the 1-cocycle relation
$\Pi(aa')=a*\Pi(a')+\Pi(a)$.\footnote{Note that this definition of a
1-cocycle 
differs from the one used in \cite{ESS} by the transformation 
$\Pi(a)\to \Pi(a^{-1})$.} 

Let $\varepsilon:X\to \mathbb C[X]^*$ be the evaluation map. 
Restriction of its values to $V$
gives a map $\tilde\varepsilon:X\to V^*$.

For $x,y\in X$  we define 
$x\circ y = \rho(\Pi^{-1}(-\tilde\varepsilon(x)))y$
(the right hand side is a map of the formal disk to $X$).  
Define\ a \ homomorphism \ of\ algebras\qquad
$R:\mathbb C[X\times X]\to \mathbb C[X\times X][[\hbar]]$ by the formula 
\begin{equation}
\label{solution}
(RF)(x,y)=F((y\circ)^{-1} x, ((y\circ)^{-1} x)\circ y),
\end{equation}
where $(y\circ)^{-1}$ denotes the inverse operator to the action of 
$y$ by $\circ$.
It is easy to see that $R=1+O(\hbar)$.

Now we will use the following result from \cite{ESS}, Section 2.4.

{\bf Proposition.} Let $G$ be a group acting on a set $X$ {\it on the
right}
via a map $\rho:G\to Aut(X)$.
Let 
$A$ be an abelian group with {\it a left} $G$-action, 
and $\pi:G\to A$ a bijective 1-cocycle. Let $\phi:X\to A$ be a
$G$-antiinvariant 
map, i.e for any $g\in G,\ x\in X$ $g\phi(x)=\phi(\rho(g^{-1})x)$. Define
$R:X\times X\to X\times X$
by 
\begin{equation}
\label{solution}
R(x,y)=((y\circ)^{-1} x, ((y\circ)^{-1} x)\circ y),
\end{equation}
where $x\circ y=\rho(\pi^{-1}(\phi(x)))y$.
Then $R$
satisfies the unitarity condition and the quantum
Yang-Baxter equation. 

{\bf Remark.} \cite{ESS} dealt with a left action of the group $G$ on X.
The above proposition is a reformulation of the corresponding result 
from Section 2.4 of \cite{ESS} in terms of the right $G$-action on X. 

This proposition (or, more precisely, its version for formal groups)
implies that $R$ is a geometric quantum R-matrix. 

It is easy to compute directly that $R=1+\hbar r+O(\hbar^2)$. 
Thus, $R$ is a quantization of $r$. The theorem is proved. $\square$

\section{Example} 

Let us show that for Examples 1 and 2, the procedure of the previous section 
gives the same quantizations as in Examples 3 and 4. 

For Example 1, this is clear: the Lie algebra $\g_+$ is 1-dimensional, and 
the computation is trivial. So let us consider Example 2. 

We have: $X$ is a finite dimensional algebra, and
 $r_c(x,y)=(xcy,-ycx)$. In this case it is easy to check 
that $\g_+$ is the right ideal $cX$ generated by $c$, with commutator 
given by $[a,b]=ab-ba$.  
The representation $V^*$ of $\g_+$ is $\g_+$ itself, with 
$a*b=-ba$. The bijective 1-cocycle $\pi$ has the form $\pi(a)=a$. 
Let us exponentiate the cocycle $\pi$. We have 
$\bar\pi(a)=(a,a)$, and 
$$
\bar\Pi(e^{\hbar a})=e^{\hbar(a,a)}=(e^{\hbar a},\frac{e^{\hbar a}-1}{a}*a)=
(e^{\hbar a},1-e^{-\hbar a}).
$$
Thus, 
$\Pi(A)=\hbar^{-1}(1-A^{-1})$. 

The map $\tilde\varepsilon$ has the form $\tilde\varepsilon(x)=cx$. Thus, 
$x\circ y=\rho(\Pi^{-1}(-cx))y=\rho((1+\hbar cx)^{-1})y=y(1+\hbar
cx)^{-1}$.
Therefore, we get 
$$
(RF)(x,y)=F(x(1+\hbar cy),y(1+\hbar cx+\hbar^2cxcy)^{-1}),
$$ 
which coincides with the geometric quantum R-matrix of Example 4. 

\section{Geometric classical r-matrices on the line and their quantization}

\begin{theorem}
Let $r$ be a geometric classical r-matrix on the affine line, which is not 
a permutation r-matrix. Then $r$ reduces to 
$r(n)=xy^n\frac{\dd}{\dd x}-yx^n\frac{\dd}{\dd y}$ for some 
$n\ge 1$, after a linear 
change of variables.  
\end{theorem}

{\it Proof}\quad  Let $\g_+$ be the Lie algebra in $\C[x]\frac{\dd}{\dd
x}$ 
associated to $r$. The $\g_+$-module $\C[x]$ has a 
nonzero finite dimensional submodule $V$. If $r$ is not 
a permutation matrix, this submodule cannot consist only of constants.
This implies that $\g_+$ is a Lie subalgebra 
of the 2-dimensional Lie algebra spanned by $\dd/\dd x$ and 
$x\dd/\dd x$ (if $\g_+$ has an element $(x^i+..)\dd/\dd x$ with $i>1$ 
then it generates an infinite dimensional space from any nonconstant 
polynomial). Consider two cases. 

1. $\g_+$ is 2-dimensional. Then $V=<1,x>$, and one can check 
that $V^*$ is isomorphic to the adjoint representation of $\g_+$. 
So if $r$ with such $\g_+$ exists, then there must exist a nondegenerate 
1-cocycle from $\g_+$ to $\g_+$, i.e. a nondegenerate derivation of $\g_+$. 
It is easy to check that such a derivation does not exist, so this case is 
impossible. 

2. $\g_+$ is 1-dimensional. Then it is 
 spanned by an element $(ax+b)\dd/\dd x$, 
such that $a\ne 0$. After a change of variable we can assume 
that $\g_+$ is spanned by $x\dd/\dd x$. Then $V=\g_-$ 
has to be spanned by $x^n$ for some $n\ge 1$, so 
$r=c(xy^n\frac{\dd}{\dd x}-yx^n\frac{\dd}{\dd y})$.

This proves the proposition, as $c$ can be easily scaled out. 
$\square$

Now let us discuss quantization of $r(n)$.
If $n=1$, this was done above: $r(1)=r_c$ from Example 2, 
where $X=\C$, and $c=1$. So the quantization 
given by the procedure of Section 2 is 
$R(1)(x,y)=(x(1+\hbar y),\frac{y}{1+\hbar x+\hbar^2 xy})$. 

If $n>1$, let us make a change of variables $x\to x^n,y\to y^n$. 
This change maps $r(n)$ to $r(1)/n$. Applying the same change of variables 
to the quantum R-matrix, we get the following quantization of $r(n)$. 
$$
R(n)(x,y)=(x(1+n\hbar y^n)^{1/n},y(1+n\hbar x^n+n^2\hbar^2x^ny^n)^{-1/n}).  
$$
It is easy to check that the procedure of Section 2 gives the same answer. 
$$ $$

{\bf \Large \noindent Acknowledgments}
$$ $$
We are thankful to the anonymous referee for useful suggestions and to 
Travis Schedler for fruitful discussions.

\end{document}